\documentclass{amsart}

\calclayout

\usepackage{amssymb, amsmath}
\usepackage{mathrsfs}
\usepackage{mathtools}
\usepackage{amscd}
\usepackage{verbatim}
\usepackage{soul}

\usepackage{wasysym}

\usepackage[utf8]{inputenc}
\usepackage{enumerate}
\usepackage{url}

\usepackage[colorinlistoftodos,prependcaption,textsize=tiny]{todonotes}

\usepackage[colorlinks,linkcolor={blue},citecolor={blue},urlcolor={red},]{hyperref}

\allowdisplaybreaks

\theoremstyle{plain}
\newtheorem{theorem}{Theorem}
\theoremstyle{remark}
\newtheorem{remark}[theorem]{Remark}

\theoremstyle{plain}

\newtheorem{lemma}[theorem]{Lemma}

\newtheorem{assumption}[theorem]{Assumption}

\numberwithin{theorem}{section}

\numberwithin{equation}{section}

\def\N{{\mathbb N}}

\def\R{{\mathbb R}}

\newcommand{\E}{{\mathbb E}}
\renewcommand{\P}{{\mathbb P}}

\newcommand{\Var}{\mathrm{Var}}

\newcommand{\mbf}[1]{\mathbf{#1}}

\begin{document}

\title[Total progeny asymptotics in multi-type branching processes]{A limit theorem for the total progeny distribution of multi-type branching processes}

\author{Jochem Hoogendijk}
\address{Mathematical institute\\
    Utrecht University\\ 3508 TA Utrecht\\ The
    Netherlands}\email{j.p.c.hoogendijk@uu.nl}

\author{Ivan Kryven}
\address{Mathematical institute\\
    Utrecht University\\ 3508 TA Utrecht\\ The
    Netherlands}\email{i.v.kryven@uu.nl}

\author{Rik Versendaal}
\address{Delft Institute of Applied Mathematics\\
    Delft University of Technology\\ 2600 GA Delft}\email{r.versendaal@tudelft.nl}

\begin{abstract}
A multi-type branching process is defined as a random tree with labeled vertices, where each vertex produces offspring independently according to the same multivariate probability distribution. We demonstrate that in realizations of the multi-type branching process, the relative frequencies of the different types in the whole tree converge to a fixed ratio, while the probability distribution for the total size of the process decays exponentially. The results hold under the assumption that all moments of the offspring distributions exist. The proof uses a combination of the arborescent Lagrange inversion formula, a measure tilting argument, and a local limit theorem. We illustrate our concentration result by showing applications to random graphs and multi-component coagulation processes.
\\
\end{abstract}

\maketitle

\section{Introduction}
A branching process is a rooted tree in which the number of offspring for each node are independent copies of a given non-negative discrete random variable, known as the offspring distribution. Originally introduced to model genealogical trees, branching processes have become central objects in probability theory due to the rich behavior exhibited by their properties. For example, the  number of nodes at a given distance from the root, the generation, may either vanish (extinction) or approach infinity (non-extinction) depending on the offspring distribution \cite{harris_1963}. The total number of nodes in a realization of the process, the total progeny, gives rise to interesting implicit equations in terms of generating functions \cite{good_1955,good_1960,good_1975,dwass_1969}. 
Branching processes are also known for their unexpected appearances in the analysis of other mathematical objects. For example, in spanning trees of random graphs \cite{bordenave_2015}, coagulation processes \cite{aldous1999,hoogendijk_2024}, and even partial differential equations \cite{mckean_1975, hoogendijk_2023}. Moreover, the simple branching process has given rise to a plethora of different generalizations, see \cite{jagers_2015} and references therein.

Our focus is on the multi-type branching process, which is a generalization in which the nodes of the tree are labeled with types. In particular, the numbers and types of offspring of a node will depend on its own type. For such processes, no explicit expression for the distribution of the total progeny is known. However, some information can be inferred from implicit formulas in terms of generating functions \cite{good_1975} or random walks  \cite{chaumont_2016}.
In this paper, we leverage on the later two results to study the
 relative frequencies of types in the total progeny and give an explicit characterization of the limiting behavior of this quantity.


Asymptotic analysis of multi-type branching processes \cite{yakovlev2009relative,Yakovlev2510} reveals that the relative frequencies of distinct types at a given generation are asymptotically multivariate normal when the initial number of ancestors is
large.
Asymptotic normality is also preserved in some modifications of the standard multi-type branching process
\cite{KOLESKO2023} and also when the tree is conditioned on non-extinction \cite{JANSON2004}.
In this paper, we show that when the relative frequencies are counted in the whole total progeny, that is, in the entire tree instead of in a fixed generation, the concentration has a more involved form. The main distinction from the former case is that, once a particular instance of the process is fixed, the offspring numbers exhibit dependencies that affect the concentration result even in the limit of diverging total number of vertices.

 Our main result gives the  following concentration of the relative frequencies of types. Let $m \in \mathbb{N}$ be the number of types and let $\mathbf{X}_k$ be the offspring random vector for a node of type $k\in[m]$. Furthermore, let the vector $\mathbf{T}$ denote the total progeny of the multi-type branching process, with the elements of this vector counting distinct types. We show that 
\begin{equation}\label{eq:eq}
    \lim\limits_{N \to\infty} \frac{1}{|\mbf{n}_N|} \log \P(\mathbf{T} = \mbf{n}_N) = -\Gamma(\boldsymbol \rho), 
\end{equation}
where 
$\Gamma(\boldsymbol \rho) =\sup_{\boldsymbol{\lambda} \in \R^m} \{\boldsymbol{\lambda} \cdot \boldsymbol{\rho} - \sum_{k=1}^m \rho_k \log \mathbb{E}[e^{\boldsymbol{\lambda} \cdot \mathbf{X}_k}]\}$. 
Here $(\mbf{n}_N)_{N \geq 1}$ with $\mbf{n}_N \in \N_0^m$ is a convergent sequence, with $\boldsymbol{\rho} = \lim_{N \to \infty} \mbf{n}_N/|\mbf{n}_N|$.  Informally, this means that depending on an arbitrarily chosen test distribution of relative frequencies $\boldsymbol{\rho}$, the probability $\mathbb{P}(\mbf{T} = \mbf{n})$ decays exponentially in $|\mbf{n}|$ with rate  $-\Gamma(\boldsymbol{\rho})$. 
In other words, the relative frequencies of types in large progenies decay exponentially while concentrating on a preferred direction $\boldsymbol \rho^*=\arg\max_{\rho\in \Delta_m} \left(-\Gamma(\boldsymbol\rho)\right)$.
Rigorously formalized in Theorem \ref{thm:main_theorem}, the result holds under certain technical assumptions, the most stringent of which is the existence of all moments of the offspring distribution.

The proof of \eqref{eq:eq} relies on finding an appropriate upper and lower bound for $\log \P(\mathbf{T} = \mbf{n}_N)$. The upper bound is found using standard Lagrange inversion formulas to relate $\P(\mathbf{T} = \mathbf{n}_N)$ to sums of random variables which can be controlled on an exponential scale using Chernoff's inequality. The lower bound is more involved, and standard Lagrange inversion techniques fail. Instead, we use arborescent Lagrange inversion \cite{bender_1998, bousquet_2003} to connect to a sum of random variables, the asymptotics of which are carefully analysed using a tilting argument and a variant of a local limit theorem \cite{petrov_1975}.

 Our main result has interesting implications for other probabilistic models. In the \emph{inhomogeneous} Erdős-Rényi random graph \cite{van_der_hofstad_2024}, the local limit is given by a unimodular multi-type branching process with Poisson offspring distribution \cite[Remark 2.13]{bollobas_2007}. Therefore, our main result characterises the asymptotic properties of the size distribution of connected components in the subcritical inhomogeneous Erdős–Rényi random graph. The large deviation principle for connected components in this random graph was also studied in \cite{andreis_2023}.
The \emph{directed colored configuration model},  features a local weak limit which can be identified with a multi-type branching process \cite[Section 4.4]{bordenave_2015}. Therefore, our result characterizes the frequency of colors in large components in this model.
The multi-type \emph{coagulation processes} are known to exhibit \emph{localization} -- a phenomenon in which the composition of large particles concentrates at a specific ratio profile \cite{ferreira2024,ferreira2021}. When the coagulation kernel is bilinearly multiplicative, localization is a special case of our result with a multi-type Poisson offspring distribution, which we recently explored in \cite[Corollary 2.8]{hoogendijk_2024}. Coagulation processes are fundamental across scientific disciplines, with applications spanning chemistry, biology, and astrophysics \cite{aldous1999}. Our findings demonstrate that multi-component clusters, or equivalently the coupled realizations of multi-type branching processes, commonly encountered in these fields, often feature a limiting composition profile that we can explicitly characterize. Previously, the analysis of such profiles relied exclusively on specially designed numerical techniques \cite{elm2020}.

In the next section, we introduce the notation and terminology, after which the main result is stated. The proof of the main result is given in Section \ref{sec:proof}.




\section{Notation and main result} \label{sec:main}
Let $m \in \N$ and consider the set of types, $[m] = \{1, \ldots, m\}$. 
We denote the offspring distribution of an individual of type $i$ by $\mathbf{X}_i = (X_{i, 1}, \ldots, X_{i, m})$, where $X_{i, j}$ with values in $\N_0$ is the number of offspring of type $j \in [m]$.  Here $\mathbf{X}_i$ are independently distributed random variables for each $i \in [m]$.  We also use the notation $|\mathbf{n}| = n_1 + \ldots + n_m$ for $\mbf{n} \in \N_0^m$. The probability generating function of $\mbf{X}_i$ is defined as $$G_{\mbf{X}_i}(\mbf{s}) := \sum\limits_{\mbf{n} \in \N_0^m} \P(\mbf{X}_i = \mbf{n}) \mbf{s}^\mbf{n},$$ where $\mbf{s}^{\mbf{n}} := s_1^{n_1} \cdot \ldots \cdot s_m^{n_m}$ for any $\mbf{n} \in \N_0^m$. To extract coefficients from the generating function, we use the notation $[\mbf{s}^{\mbf{n}}] G_{\mbf{X}}(\mbf{s}) := \P(\mbf{X} = \mbf{n})$.


 The multi-type branching process is the sequence $\{\mbf{Z}_n\}_{n\in\mathbb{N}}$ defined iteratively by the recurrence equation. Choose the type $i \in [m]$ of the root and set  $\mbf{Z}_0 = \mbf{e}_i$ to be the unit vector. Then
 $$\mbf{Z}_{n+1} = \sum\limits_{k=1}^m \sum\limits_{j=1}^{Z_{n, k}} \mbf{X}_{k, (n, j)},$$ where $\mbf{X}_{k, (n, j)} \sim \mbf{X}_k$ for all $n, j \in \N$.
 The total progeny of the process is defined as $\mbf{T}^{(i)} = \sum_{n=0}^\infty \mbf{Z}_n.$ Suppose we choose a root of random type,  $\P(\mathbf{Z}_0 = \mbf{e}_i) = p_i$, according to some probability distribution $\{p_i\}_{i \in [m]}$. We then define the total progeny as $\mbf{T} = \sum_{n=0}^\infty \mbf{Z}_n$. It follows that $\P(\mbf{T} = \mbf{n}) = \sum_{i=1}^m p_i \P(\mbf{T}^{(i)} = \mbf{n})$. Finally, we use the notation $\Delta_m$ for the $m-1$ dimensional simplex $\Delta_m = \{\boldsymbol{\rho} \in \R^m : \rho_1 + \ldots + \rho_m  = 1, \boldsymbol{\rho} \geq 0\}$. 

To state the main theorem, we need the following assumptions.
\begin{assumption}\label{ass:sequence_assumption}
    Let $(\mbf{n}_N)_{N \geq 1}$ be a sequence in $\N_0^m$ such that there exists a constant $C \geq 0$ and a stochastic vector $\boldsymbol{\rho} \in \R^m$ such that $|\mbf{n}_N| \to \infty$ as $N \to \infty$ and $|\mathbf{n}_N - |\mbf{n}_N| \boldsymbol{\rho}| \leq C \sqrt{|\mathbf{n}_N|}.$ 
\end{assumption}
\begin{assumption}\label{ass:offspring_assumptions}
    The collection of offspring distributions $\{\mbf{X}\}_{k \in [m]}$ satisfies:
    \begin{enumerate}[a)]
        \item There exists $r > 0$ such that $\E[e^{\boldsymbol{\lambda} \cdot \mbf{X}_k}] < \infty$ for $\boldsymbol{\lambda} \in B_r(\mathbf{0})$ for each $k \in [m]$.
        \item For all $i, j \in [m]$ and all $n \in \N_0$, $\P(X_{i,j} = n) > 0$, {\it i.e.} $X_{i,j}$ are fully supported. 
    \end{enumerate}
\end{assumption}
\begin{theorem}\label{thm:main_theorem}
    Consider a multi-type branching process as described above with the offspring satisfying Assumption \ref{ass:offspring_assumptions}. Let $(\mbf{n}_N)_N$ be a sequence in $\N_0^m$ satisfying Assumption \ref{ass:sequence_assumption}. Then,
    \begin{equation}
        \lim\limits_{N \to \infty} \frac{1}{|\mbf{n}_N|} \log \P(\mbf{T} = \mbf{n}_N) = - \Gamma(\boldsymbol{\rho}),
    \end{equation}
    where $\Gamma: \Delta_m \to \R$ with
    \begin{equation}\label{eq:GammaRho}
        \Gamma(\boldsymbol{\rho}) = \sup\limits_{\boldsymbol{\lambda} \in \R^m} \{\boldsymbol{\lambda} \cdot \boldsymbol{\rho} - \sum\limits_{k=1}^m \rho_k \log \E[e^{\boldsymbol{\lambda} \cdot \boldsymbol{X}_k}]\}.
    \end{equation}
\end{theorem}

 Let $\mathbf{T}_N := \mathbf{T} \mid |\mathbf{T}| = N$ be the total progeny vector conditioned on having $N$ individuals. One implication of the above theorem is that $\mathbf{T}_N/N$ converges in distribution to $\boldsymbol{\rho}^* = \mathrm{argmin}_{\boldsymbol{\rho} \in \Delta_m} \Gamma(\boldsymbol{\rho})$. In general, it is not possible to obtain an explicit expression for $\boldsymbol{\rho}^*$. However, if the average offspring matrix $\mathbb{E} (\mbf{X}_1, \ldots, \mbf{X}_m)^T$ is a right stochastic matrix, it can be shown that $\boldsymbol{\rho}^*$ is the principal eigenvector of this matrix.

We also note that Theorem \ref{thm:main_theorem} seems reminiscent of a large deviation principle (LDP). However, it is not generally the case for an LDP of a sequence of random variables $(X_n)_{n\geq 1}$, that there are sets $A$ for which $(1/n) \log \P(X_n \in A)$ converges. The statement we prove is about the singleton sets, which can not be extracted from an LDP setting because singletons have empty interior. Notice that a less stringent statement about the quantities $\P(\mbf{T} \geq \mbf{n})$ and $\P(\mbf{T} \leq \mbf{n})$ can be obtained by modifying part of our argument and applying the G\"{a}rtner-Ellis theorem \cite{dembo_1998}. In this case, $\Gamma(\boldsymbol{\rho})$ as defined in equation \ref{eq:GammaRho}, indeed acts as the rate function.

We now give some comments on the assumptions made for the theorem. First of all, note that that the assumption $\P(X_{i,j} = 0) > 0$, contained in Assumption \ref{ass:sequence_assumption} has to be made, because otherwise $\mbf{T}$ is almost surely infinite, {\it i.e.} $\P(|\mbf{T}| = \infty) = 1$. Furthermore, Assumption \ref{ass:offspring_assumptions}b) is not strictly necessary, but rather made to simplify the exposition. It will allow us to apply the multivariate local limit theorem from \cite[Theorem 7.2]{petrov_1975} in the proof. If one is interested in applying the theorem under less stringent conditions on the support of the offspring distribution, the proof we give can be changed in combination with the local limit theorem.  We also believe that the assumption $|\mbf{n}_N - |\mbf{n}_N| \boldsymbol{\rho}| \leq C \sqrt{|\mbf{n}_N|}$ (see Assumption \ref{ass:sequence_assumption}) can be improved if one has a stronger or more explicit version of a multivariate local limit theorem.

Finally, we provide some intuition for Theorem \ref{thm:main_theorem}.
In the subcritical case, one my view the tail of total progeny distribution as a conditional survival probability. Hence, it is expected that this tail is exponential, the rate of which Theorem \ref{thm:main_theorem} quantifies precisely. 
However, it seems surprising that Theorem \ref{thm:main_theorem} should also hold in the supercritical setting. The reason is that in the supercritical regime, a `large' tree is simply more likely to keep producing individuals instead of dying out. Therefore, the probability of seeing a \emph{finite} large tree still decays exponentially.

\section{Proof of the main result}\label{sec:proof}
The proof consists of two parts. We prove an asymptotic upper bound and lower bound for $(1/|\mbf{n}_N|)\log \P(\mbf{T} = \mbf{n}_N)$, the gap between which becomes tight when $N$ tends to infinity. To obtain an upper bound, we  use Lagrange-Good inversion (Lemma \ref{lem:good_lagrange}), interpret the resulting quantity probabilistically and use Chernoff's inequality. The lower bound  involves the arborescent Lagrange inversion formula (Lemma \ref{lem:arborescent_good_lagrange}), which allow us to circumvent the determinant term that is present in the standard Lagrange-Good inversion and is otherwise hard to control. After tilting the probability measure in the right direction, we conclude by applying a multivariate local limit theorem. This strategy parallels and extends the standard structure of proofs in classical large deviations theory (see for example \cite{den_hollander_2000,dembo_1998}), where the upper bound uses Chernoff's inequality and the lower bound is obtained by a tilting procedure.
For notational convenience, we will mostly write $\mbf{n}$ instead of $\mbf{n}_N$ in the sequel, though the latter is always implied. 

\subsection{Lagrange inversion}
We start with stating several auxiliary  results. 
\begin{lemma}[\cite{good_1955}]\label{lem:implicit_system}
    Consider a multi-type branching processes as described in the previous section. Then, the following system of equations holds:
    \begin{equation}
        \begin{split}
            G_{\mbf{T}_1}(\mbf{s}) &= s_1 G_{\mbf{X}_1}(G_{\mbf{T}_1}, \ldots, G_{\mbf{T}_m}),\\
            &\stackrel{\vdots}{\hphantom{=}}\\
            G_{\mbf{T}_m}(\mbf{s}) &= s_m G_{\mbf{X}_m}(G_{\mbf{T}_1}, \ldots, G_{\mbf{T}_m}),
        \end{split}
    \end{equation}
    in the formal power series sense.
\end{lemma}

Lemma \ref{lem:implicit_system} can be conveniently combined with Lagrange inversion to obtain an expression for $\P(\mbf{T} = \mbf{n})$. We will use two variants of Lagrange inversion, which we now state below. The first form will be useful for the upper bound, whereas the second form will be used in the proof of the lower bound.

\begin{lemma}[Lagrange-Good inversion, \cite{good_1960, bergeron_1998}]\label{lem:good_lagrange}
    For any formal power series $F(\mathbf{s})$ in the variables $(s_1, \ldots, s_m)$ and for all $n \in \mathbb{N}_0^m$,
    \begin{equation}
        [\mathbf{s}^\mathbf{n}]F(G_{\mbf{T}_1}(\mbf{s}), \ldots, G_{\mbf{T}_{m}}(\mbf{s})) = [\mathbf{r}^\mathbf{n}]F(\mbf{r}) \det(K(\mathbf{r})) G_{\mathbf{X}_1}^{n_1} \cdot \ldots \cdot G_{\mathbf{X}_m}^{n_m}
    \end{equation}
    where
    \begin{equation}\label{eq:lagrange_inversion_multi_type}
        K(\mbf{r}) := \left[\delta_{i,j} - \frac{r_i}{G_{\mbf{X}_i}(r_1, \ldots, r_m)} \frac{\partial G_{\mbf{X}_i}}{\partial r_j}\right]_{1 \leq i, j \leq m},
    \end{equation}
    is a matrix associated to the generating functions $G_{\mbf{X}_1}, \ldots, G_{\mbf{X}_m}$.
\end{lemma}
\begin{remark}
    On choosing $F(s_1, \ldots, s_m) = \sum_{k=1}^m p_k s_k$ in Lemma \ref{lem:good_lagrange}, we obtain the formula
    \begin{equation}
        \P(\mbf{T} = \mbf{n}) = \sum\limits_{k=1}^m p_k [\mbf{r}^{\mbf{n}}] r_k \mathrm{det}(K(\mbf{r})) G_{\mbf{X}_1}^{n_1} \cdot \ldots \cdot G_{\mbf{X}_m}^{n_m}.
    \end{equation}
\end{remark}
Another form of Lagrange inversion that we will rely on is the arborescent Lagrange inversion formula \cite{bender_1998, bousquet_2003}. In order to formulate it, we have to introduce the notion of tree derivative \cite{bender_1998}. Let $\mathcal{D}$ be a directed graph with vertex set $V$ and edge set $E$ and let the vectors $\mbf{x}$ and $\mbf{f}(\mbf{x})$ be indexed by $V$. The derivative with respect to $\mathcal{D}$ is
\begin{equation}
    \frac{\partial \mbf{f}}{\partial \mathcal{D}} := \prod\limits_{j\in V}\left\{\left(\prod\limits_{(i, j) \in E} \frac{\partial}{\partial x_i}\right) f_j(\mbf{x})\right\}.
\end{equation}
Then,
\begin{lemma}[Arborescent Lagrange inversion, \cite{bender_1998, bousquet_2003}]\label{lem:arborescent_good_lagrange}
    For any formal power series $F(\mathbf{s})$ in the variables $(s_1, \ldots, s_m)$ and for all $\mbf{n} \in \mathbb{N}_0^m$,
    \begin{equation}
        [\mathbf{s}^\mathbf{n}]F(G_{\mbf{T}_1}, \ldots, G_{\mbf{T}_{m}}) = \frac{1}{\prod_{i=1}^m n_i} [\mbf{r}^{\mbf{n}-\mbf{1}}] \sum\limits_{\mathcal{T}} \frac{\partial(F, G_{\mbf{X}_1}^{n_1}, \ldots, G_{\mbf{X}_m}^{n_m})}{\partial \mathcal{T}},
    \end{equation}
    where the sum is taken over all trees $\mathcal{T}$ with $V = \{0, 1, \ldots, m\}$ and edges directed towards $0$.
\end{lemma}
The particularly convenient feature of the arborescent Lagrange inversion formula is that all individual terms in the sum are non-negative when the coefficients of the involved generating functions are non-negative. This is not the case in Lemma \ref{lem:good_lagrange}, where the determinant can give negative terms in its expansion. The non-negativity of the coefficient in arborescent Lagrange inversion will turn out to be a crucial property for proving the lower bound. 

\subsection{Proof of the upper bound}
We will show in this section that
\begin{equation}
    \limsup\limits_{N \to \infty} \frac{1}{|\mathbf{n}_N|} \log \P(\mathbf{T} = \mbf{n}) \leq -\Gamma(\boldsymbol{\rho}).
\end{equation}
\begin{proof}
 As as consequence of Lemma \ref{lem:good_lagrange}, we have the formula
\begin{equation}\label{eq:total_progeny_1}
    \P(\mbf{T} = \mbf{n}) = \sum\limits_{k=1}^m p_k [\mbf{r}^{\mbf{n}}] r_k \mathrm{det}(K(\mbf{r})) G_{\mbf{X}_1}^{n_1} \cdot \ldots \cdot G_{\mbf{X}_m}^{n_m}.
\end{equation}
Before analyzing this inequality further, we will bound the determinant term. First, notice that for any $i, j \in [m]$ and $\mbf{r} > \mbf{0}$,
\begin{equation}\label{eq:determinant_term_est_1}
    |K(\mbf{r})_{i, j}| = \left|\delta_{i,j} - \frac{r_i}{G_{\mbf{X}_i}}\frac{\partial G_{\mbf{X}_i}}{\partial r_j}\right| \leq 1 + \frac{r_i}{G_{\mbf{X}_i}} \frac{\partial G_{\mbf{X}_i}}{\partial r_j}
\end{equation}
since $$\frac{r_i}{\partial G_{\mbf{X}_i}} \frac{G_{\mbf{X}_i}}{\partial r_j} > 0$$ for $\mbf{r} > \mbf{0}$. Using estimate \eqref{eq:determinant_term_est_1}, and writing $S_m$ for the symmetric group on $[m]$, we see that 
\begin{equation}\label{eq:determinant_est_1}
    \begin{split}
    \det (K(\mbf{r})) &= \sum\limits_{\sigma \in S_m} \mathrm{sgn}(\sigma) \prod\limits_{i=1}^m (K(\mbf{r}))_{i, \sigma(i)}\\
    &\leq \sum\limits_{\sigma \in S_m} \prod\limits_{i=1}^m |(K(\mbf{r}))_{i, \sigma(i)}|\\
    &\leq \sum\limits_{\sigma \in S_m} \prod\limits_{i=1}^m \left(1 + \frac{r_i}{G_{\mbf{X}_i}} \frac{\partial G_{\mbf{X}_i}}{\partial r_{\sigma(i)}}\right)
    \end{split}
\end{equation}
where we used the triangle inequality in  the second and third steps. Consider the following formulation of `expanding the brackets':
\begin{lemma}\label{lem:bracket_expansion}
    Consider some finite set $\mathcal{V}$ and let $f : \mathcal{V} \to \R$. Then, 
    \begin{equation}
        \prod\limits_{v \in \mathcal{V}} (1 + f(v)) = \sum\limits_{V \subseteq \mathcal{V}} \prod\limits_{v \in V} f(v).
    \end{equation}
\end{lemma}
\noindent Applying Lemma \ref{lem:bracket_expansion} to \eqref{eq:determinant_est_1}
results in the estimate
\begin{equation}\label{eq:determinant_est_2}
    \mathrm{det}(K(\mbf{r})) \leq \sum\limits_{\sigma \in S_m} \sum\limits_{J \subseteq [m]} \prod\limits_{i \in J} \frac{r_i}{G_{\mbf{X}_i}} \frac{\partial G_{\mbf{X}_i}}{\partial r_{\sigma(i)}}.
\end{equation}
Combining \eqref{eq:total_progeny_1} and \eqref{eq:determinant_est_2} yields for all $\mbf{n}$ such that $|\mbf{n}| > M$ with $M$ large enough that that
\begin{equation}\label{eq:intermediate_step_1}
    \P(\mbf{T} = \mbf{n}) \leq \sum\limits_{k=1}^m \sum\limits_{\sigma \in S_m} \sum\limits_{J \subseteq [m]} p_k [\mbf{r}^{\mbf{n}}] r_k \prod\limits_{i \in J} \frac{r_i}{G_{\mbf{X}_i}} \frac{\partial G_{\mbf{X}_i}}{\partial r_{\sigma(i)}} G_{\mbf{X}_1}^{n_1} \cdot \ldots \cdot G_{\mbf{X}_m}^{n_m}.
\end{equation}
We introduce the random variable $\mbf{Y}_{i, j}$ with probability generating function
\begin{equation}
    G_{\mbf{Y}_{i, j}} := \frac{1}{C_{i, j}} \frac{\partial G_{\mbf{X}_i}}{\partial r_j},
\end{equation}
where $C_{i, j}$ is a normalization constant. Now, 
\begin{equation}
    \begin{split}
        &[\mbf{r}^{\mbf{n}}]r_k \prod\limits_{i \in J} \frac{r_i}{G_{\mbf{X}_i}} \frac{\partial G_{\mbf{X}_i}}{\partial r_{\sigma(i)}} G_{\mbf{X}_1}^{n_1} \cdot \ldots \cdot G_{\mbf{X}_m}^{n_m}\\
        &= \P\left(\sum\limits_{i_1 = 1}^{n_1 - \mathbf{1}_J(1)} \mbf{X}_1^{(i_m)}+ \ldots + \sum\limits_{i_m = 1}^{n_m - \mbf{1}_J(m)} \mbf{X}_m^{(i_m)} + \sum\limits_{i \in J} \mbf{Y}_{i, \sigma(i)} = \mbf{n} - \mbf{e}_k - \sum\limits_{i \in J} \mbf{e_i}\right),
    \end{split}
\end{equation}
with $\mathbf{1}_J$ being the indicator function of the set $J$. The terms $\mbf{X}_1^{(i_m)}, \ldots \mbf{X}_m^{(i_m)}$ are i.i.d. copies of $\mbf{X}_1, \ldots, \mbf{X}_m$. The above equality follows by the fact that a product of probability generating function becomes a probability generating function of the sum of the associated random variables, and the variables $r_k$ and $r_i$ induce a shift, resulting in the term $\mathbf{e}_k + \sum_{i \in J} \mathbf{e}_i$. Therefore, we may interpret the inequality \eqref{eq:intermediate_step_1} as
\begin{equation}
    \begin{split}
        \P(\mbf{T} = \mbf{n}) &\leq \sum\limits_{k=1}^m \sum\limits_{\sigma \in S_m} \sum\limits_{J \subseteq [m]} p_k \left(\prod\limits_{i \in J} \frac{1}{C_{i, \sigma(i)}}\right)\\
        &\quad \cdot \P\Bigg(\sum\limits_{i_1 = 1}^{n_1 - \mathbf{1}_J(1)} \mbf{X}_1^{(i_m)}+ \ldots + \sum\limits_{i_m = 1}^{n_m - \mbf{1}_J(m)} \mbf{X}_m^{(i_m)} + \sum\limits_{i \in J} \mbf{Y}_{i, \sigma(i)} = \mbf{n} - \mbf{e}_k - \sum\limits_{i \in J} \mbf{e_i}\Bigg),
    \end{split}
\end{equation}
Applying Chernoff's inequality and absorbing the terms not depending on $\mathbf{n}$ into the constant $C$, we obtain 
\begin{equation}
    \P(\mbf{T} = \mbf{n}) \leq C \exp\left(-\boldsymbol{\lambda} \cdot \mbf{n}\right) \E[e^{\boldsymbol{\lambda} \cdot \mbf{X}_1}]^{n_1} \cdot \ldots \cdot \E[e^{\boldsymbol{\lambda} \cdot \mbf{X}_m}]^{n_m}.
\end{equation}
Taking logarithms on both sides and dividing by $|\mbf{n}|$ yields
\begin{equation}
    \frac{1}{|\mbf{n}|} \log \P(\mbf{T} = \mbf{n}) \leq \frac{1}{|\mbf{n}|} \log C  -\boldsymbol{\lambda} \cdot \frac{\mbf{n}}{|\mbf{n}|} + \sum\limits_{k=1}^m \frac{{n_{k}}}{|\mbf{n}|} \log \E[e^{\boldsymbol{\lambda} \cdot \mbf{X}_k}].
\end{equation}\
Take the limit supremum on both sides over $N$ and optimizing over $\boldsymbol{\lambda}$ yields
\begin{equation}
    \limsup_{N\to \infty} \frac{1}{|\mbf{n}_N|} \log \P(\mbf{T} = \mbf{n}_N) \leq - \sup\limits_{\boldsymbol{\lambda}} \{-\boldsymbol{\lambda} \cdot \boldsymbol{\rho} + \sum\limits_{k=1}^m \rho_k \log \E[e^{\boldsymbol{\lambda} \cdot \mbf{X}_k}]\}.
\end{equation}
\end{proof}

\subsection{Proof of the lower bound}
The goal of this section is to show that
\begin{equation}
    \liminf\limits_{N \to \infty} \frac{1}{|\mathbf{n}_N|} \log \P(\mathbf{T} = \mathbf{n}_N) \geq - \Gamma(\boldsymbol{\rho}).
\end{equation}
\begin{proof}
The arborescent Lagrange inversion formula tells us that
\begin{equation}
    \P(\mbf{T} = \mbf{n}) = \frac{1}{\prod_{i=1}^m n_i}  \sum\limits_{\mathcal{T}} [\mbf{r}^{\mbf{n}-\mbf{1}}]\frac{\partial(\sum_{k=1}^m p_k r_k, G_{\mbf{X}_1}^{n_1}, \ldots, G_{\mbf{X}_m}^{n_m})}{\partial \mathcal{T}}.
\end{equation}
Note that each of the terms in the sum is non-negative, since the probability generating functions all have non-negative coefficients. Therefore, for any directed tree $\mathcal{T}$ with vertex set $V = \{0, 1, \ldots, m\}$ and edges directed towards $0$, 
\begin{equation}\label{eq:lower_bound_arbo}
    \P(\mbf{T} = \mbf{n}) \geq \frac{1}{\prod_{i=1}^m n_i} [\mbf{r}^{\mbf{n}-\mbf{1}}] \frac{\partial(\sum_{k=1}^m p_k r_k, G_{\mbf{X}_1}^{n_1}, \ldots, G_{\mbf{X}_m}^{n_m})}{\partial \mathcal{T}}.
\end{equation}
We choose the tree $\mathcal{T}$ for which there exists a directed edge from $j$ to $i$ if and only if $j - i = 1$. This graph is known as the directed path,  see Figure \ref{fig:tree} for an example for $m = 5$. 
\begin{figure}
    \centering
\begin{tikzpicture}[->, >=stealth, auto, node distance=1.5cm, thick]
    \foreach \i in {0,...,5}
        \node[circle, draw] (node\i) at (\i,0) {\i};
    
    \foreach \i [evaluate=\i as \j using {int(\i+1)}] in {0,...,4}
        \draw (node\j) -- (node\i);
\end{tikzpicture}
\caption{An example of the tree used in the proof for $m = 5$.}
\label{fig:tree}
\end{figure}
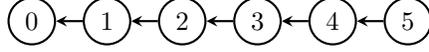

Now, we compute the tree derivative and obtain
\begin{equation}
    \begin{split}
        &\frac{\partial}{\partial \mathcal{T}} \Big(\sum_{i=1}^m p_i r_i, G_{\mbf{X}_1}^{n_1}, \ldots, G_{\mbf{X}_m}^{n_m}\Big)\\
        &= \frac{\partial}{\partial r_1}\Big(\sum\limits_{i=1}^m p_i r_i\Big) \cdot \frac{\partial}{\partial r_2}\left(G_{\mbf{X}_1}^{n_1}\right) \cdot \ldots \cdot \frac{\partial}{\partial r_{m}}\left(G_{\mbf{X}_{m-1}}^{n_{m-1}}\right) \cdot G_{\mbf{X}_m}^{n_m}\\
        &= p_1 \frac{\partial}{\partial r_2}\left(G_{\mbf{X}_1}^{n_1}\right) \cdot \ldots \cdot \frac{\partial}{\partial r_{m}}\left(G_{\mbf{X}_{m-1}}^{n_{m-1}}\right) \cdot G_{\mbf{X}_m}^{n_m}\\
        &= p_1 n_1 \left(G_{\mbf{X}_1}\right)^{n_1-1} \frac{\partial G_{\mbf{X}_1}}{\partial r_2} \cdot \ldots \cdot n_{m-1} \left(G_{\mbf{X}_{m-1}}^{n_{m-1}-1}\right) \frac{\partial G_{\mbf{X}_{m-1}}}{\partial r_{m-1}}\cdot G_{\mbf{X}_m}^{n_m}.
    \end{split}
\end{equation}
We combine the above computation with \eqref{eq:lower_bound_arbo} to arrive at
\begin{equation}\label{eq:lower_bound_arbo_2}
    \begin{split}
        \P(\mbf{T} = \mbf{n})
        &\geq \frac{p_1}{n_m} [\mbf{r}^{\mbf{n}-\mbf{1}}] \left(G_{\mbf{X}_1}\right)^{n_1-1} \frac{\partial G_{\mbf{X}_1}}{\partial r_2} \cdot \ldots \cdot \left(G_{\mbf{X}_{m-1}}^{n_{m-1}-1}\right) \frac{\partial G_{\mbf{X}_{m-1}}}{\partial r_{m}}\cdot G_{\mbf{X}_m}^{n_m}\\
        &= \frac{p_1}{n_m} [\mbf{r}^{\mbf{n}-\mbf{1}}](G_{\mbf{X}_1})^{n_1-1} \cdot \ldots \cdot (G_{\mbf{X}_m})^{n_m-1} \cdot \frac{\partial G_{\mbf{X}_1}}{\partial r_2} \cdot \ldots \frac{\partial G_{\mbf{X}_{m-1}}}{\partial r_m} \cdot G_{\mbf{X}_m}.
    \end{split}
\end{equation}
By introducing an appropriate rescaling, we can see the terms $\frac{\partial G_{\mbf{X}_1}}{\partial r_2}, \ldots, \frac{\partial G_{\mbf{X}_{m-1}}}{\partial r_m}$ as probability generating functions themselves, which we denote by $G_{\mbf{Y}_1}, \ldots G_{\mbf{Y}_{m-1}}$, respectively. We will also denote $\mbf{X}_m$ by $\mbf{Y}_m$. This means that we can write inequality \eqref{eq:lower_bound_arbo} as 
\begin{equation}
    \begin{split}
        \P(\mbf{T} = \mbf{n})
        &\geq \frac{Cp_1}{n_m}\P\left(\sum\limits_{i_1=1}^{n_1-1} \mbf{X}_1^{(i_1)} + \ldots + \sum\limits_{i_{m} = 1}^{n_{m}-1} \mbf{X}_{m}^{(i_{m})} + \mbf{Y}_1 + \ldots + \mbf{Y}_{m-1} + \mbf{Y}_m = \mbf{n}-\mbf{1}\right),
    \end{split}
\end{equation}
where $\mbf{X}_k^{(i_k)} \stackrel{d}{=} \mbf{X}_k$ and all terms are mutually independent, for each $k \in [m]$, $i_k \in \N$. Define
\begin{equation}
    P(\mbf{n}) := \P\Bigg(\sum\limits_{i_1=1}^{n_1} \mbf{X}_1^{(i_1)} + \ldots + \sum\limits_{i_{m} = 1}^{n_{m}} \mbf{X}_{m}^{(i_{m})} + \mbf{Y}_1 + \ldots + \mbf{Y}_{m-1} + \mbf{Y}_m = \mbf{n}\Bigg)
\end{equation}
and
\begin{equation}
    \mbf{X}'_k := \mbf{X}_k - \mbf{e}_k.
\end{equation}
Then, we see that
\begin{equation}
    P(\mbf{n}) = \P\Bigg(\sum\limits_{i_1=1}^{n_1} \mbf{X}'^{(i_1)}_1 + \ldots + \sum\limits_{i_{m} = 1}^{n_{m}} \mbf{X}'^{(i_{m})}_m + \mbf{Y}_1 + \ldots + \mbf{Y}_{m-1} + \mbf{Y}_m = \mbf{0}\Bigg),
\end{equation}
where $\mbf{X}'^{(i_k)}_k \stackrel{d}{=} \mbf{X}'_k$ and all terms are mutually independent, for each $k \in [m]$, $i_k \in \N$. Since the offspring components of the vectors are independent we may write 
\begin{align*}
        P(\mbf{n})
        &= \P\left(\sum\limits_{i_1=1}^{n_1} X_{1,1}^{'(i_1)} + \ldots + \sum\limits_{i_{m} = 1}^{n_{m}} X_{m, 1}^{'(i_{m})} + Y_{1,1} + \ldots + Y_{m, 1} = 0 \right) \cdot \ldots \\
        &\quad \cdot \P\left(\sum\limits_{i_1=1}^{n_1} X_{1,m}^{'(i_1)} + \ldots + \sum\limits_{i_{m} = 1}^{n_{m}} X_{m, m}^{'(i_{m})} + Y_{1,m} + \ldots + Y_{m, m} = 0 \right).
\end{align*}
Now, we work our way towards the tilting argument. Given $j \in [m]$, define
\begin{equation}\label{eq:tilting_function}
    \phi_j(\lambda):=\E[e^{\lambda X'_{1,j}}]^{\rho_1} \cdot \ldots \cdot \E[e^{\lambda X'_{mj}}]^{\rho_m}.
\end{equation}
By assumption, $X_{i,j}'$ is fully supported on $\N_0 \cup \{-1\}$, so we find that $\lim_{\lambda \to \pm \infty} \phi_j(\lambda) = \infty$. By convexity of $\phi_j$ for each $j \in [m]$, there exists a unique $\tau_j$ such that $\phi_j(\tau_j) \leq \phi_j(\lambda)$ for all $\lambda \in \R$. 

Given $i, j \in [m]$, let $F_{i,j}$ be the cumulative distribution function of $X_{i, j}'$. Then, define the tilted cumulative distribution function $\widehat{F}_{i,j}$ as
\begin{equation}\label{eq:new_distr_1}
    \widehat{F}_{i,j}(x) := \frac{1}{\E[e^{\tau_j X'_{i,j}}]} \int_{(-\infty, x]} e^{\tau_j y} {\rm d}F_{i,j}(y).
\end{equation}
Denote the associated random variable by $\widehat{X}'_{i,j}$. Similarly, we let $H_{i,j}$ be the cumulative distribution function of $Y_{i, j}$. However, we will not be tilting it. Finally, let $A(\mbf{n})$ be the set
\begin{equation}
    \begin{split}
        A(\mbf{n}) &= \Big\{x_1^{(1)}, \ldots, x_1^{(n_1)}, \ldots, x_m^{(1)}, \ldots, x_m^{(n_m)}, y_1, \ldots, y_m :\\
        &\quad \sum\limits_{i_1=1}^{n_1} x_1^{(i_1)} + \ldots + \sum\limits_{i_m = 1}^{n_m} x_m^{(i_m)} + y_1 + \ldots + y_m  = 0\Big\}.
    \end{split}
\end{equation}
Observe that for any $j \in [m]$
\begin{equation}\label{eq:lower_bound_1}
    \begin{split}
        &\P\left(\sum\limits_{i_1=1}^{n_1} X_{1,j}^{'(i_1)} + \ldots + \sum\limits_{i_{m} = 1}^{n_{m}} X_{m, j}^{'(i_{m})} + Y_{1,j} + \ldots + Y_{m, j} = 0 \right)\\
        &= \int_{A(\mbf{n})} \left(\prod\limits_{i_1=1}^{n_1}  {\rm d} F_{1j}(x_1^{(i_1)})\right) \cdot \ldots \cdot \left(\prod\limits_{i_m = 1}^{n_m}  {\rm d} F_{mj}(x_m^{(i_m)})\right)
         {\rm d} H_{1j}(y_1) \cdot ... \cdot {\rm d}H_{mj}(y_m).\\
    \end{split}
\end{equation}
Applying the tilting \eqref{eq:new_distr_1}, we see that
\begin{equation}
    \begin{split}
        &\P\left(\sum\limits_{i_1=1}^{n_1} X_{1,j}^{'(i_1)} + \ldots + \sum\limits_{i_{m} = 1}^{n_{m}} X_{m, j}^{'(i_{m})} + Y_{1,j} + \ldots + Y_{m, j} = 0 \right)\\
        &= \E[e^{\tau_1 X'_{1,j}}]^{n_1} \cdot \ldots \cdot \E[e^{\tau_m X'_{m,j}}]^{n_m}\exp\left(-\sum\limits_{k=1}^m \tau_k n_k\right)\\
        &\quad \cdot \int_{A(\mbf{n})} \left(\prod\limits_{i_1=1}^{n_1}  {\rm d} \widehat{F}_{1j}(x_1^{(i_1)})\right) \cdot \ldots \cdot \left(\prod\limits_{i_m = 1}^{n_m}  {\rm d} \widehat{F}_{mj}(x_m^{(i_m)})\right) {\rm d}H_{1j}(y_1) \cdot ... \cdot {\rm d}H_{mj}(y_m).\\
    \end{split}
\end{equation}
Interpreting the latter integral probabilistically again, we see that
\begin{equation}\label{eq:lower_bound_2}
    \begin{split}
        &\P\left(\sum\limits_{i_1=1}^{n_1} X_{1,j}^{'(i_1)} + \ldots + \sum\limits_{i_{m} = 1}^{n_{m}} X_{m, j}^{'(i_{m})} + Y_{1,j} + \ldots + Y_{m, j} = 0 \right)\\
        &= \E[e^{\tau_1 X'_{1,j}}]^{n_1} \cdot \ldots \cdot \E[e^{\tau_m X'_{m,j}}]^{n_m}\\
        &\quad \cdot \P\left(\sum\limits_{i_1=1}^{n_1} \widehat{X}_{1,j}^{'(i_1)} + \ldots + \sum\limits_{i_{m} = 1}^{n_{m}} \widehat{X}_{m, j}^{'(i_{m})} + Y_{1,j} + \ldots + Y_{m, j} = 0 \right).\\
    \end{split}
\end{equation}
We will now use the multivariate local limit theorem from \cite[Theorem 7.2]{petrov_1975} to obtain asymptotic behavior for the probability appearing on the right-hand side of the last equation. Before applying it, we have to introduce some definitions. Define 
\begin{equation}
    \widehat{P}_j(\mbf{n}) := \P\left(\sum\limits_{i_1=1}^{n_1} \widehat{X}_{1,j}^{'(i_1)} + \ldots + \sum\limits_{i_{m} = 1}^{n_{m}} \widehat{X}_{m, j}^{'(i_{m})} + Y_{1,j} + \ldots + Y_{m, j} = 0 \right)
\end{equation}
as well as
\begin{equation}
    \widehat{M}_j(\mbf{n}) := \E\left[\sum\limits_{i_1=1}^{n_1} \widehat{X}^{'(i_1)}_{1, j} + \ldots + \sum\limits_{i_m = 1}^{n_m} \widehat{X}^{'(i_m)}_{m, j} + Y_{1, j} + \ldots + Y_{m ,j}\right]
\end{equation}
and 
\begin{equation}
    \widehat{V}_j(\mbf{n}) := \Var\left[\sum\limits_{i_1=1}^{n_1} \widehat{X}^{'(i_1)}_{1, j} + \ldots + \sum\limits_{i_m = 1}^{n_m} \widehat{X}^{'(i_m)}_{m, j} + Y_{1, j} + \ldots + Y_{m ,j}\right].
\end{equation}
The latter two expressions simplify to
\begin{equation}
    \widehat{M}_j(\mbf{n}) = \sum\limits_{k=1}^m n_k \E[\widehat{X}'_{kj}] + \E[Y_{k, j}]
\end{equation}
and
\begin{equation}
    \widehat{V}_j(\mbf{n}) = \sum\limits_{k=1}^m n_k \Var[\widehat{X}'_{kj}] + \Var[Y_{k, j}].
\end{equation}
By \cite[Theorem 7.2]{petrov_1975}, we see that 
\begin{equation}
    \left|\sqrt{\widehat{V}_j(\mbf{n}_N)}\widehat{P}_j(\mbf{n}_N) - \frac{1}{\sqrt{2\pi}} \exp\left(-\frac{\widehat{M}_j(\mbf{n}_N)^2}{2 \widehat{V}_j(\mbf{n}_N)}\right)\right| \to 0
\end{equation}
as $N \to \infty$. Therefore, it follows that there exists a non-negative sequence $\epsilon_N$ with $\lim_{N\to\infty} \epsilon_N = 0$ such that 
\begin{equation}
-\epsilon_N \leq \sqrt{\widehat{V}_j(\mbf{n}_N)}\widehat{P}_j(\mbf{n}_N) - \frac{1}{\sqrt{2\pi}} \exp\left(-\frac{\widehat{M}_j(\mbf{n}_N)^2}{2 \widehat{V}_j(\mbf{n}_N)}\right) \leq \epsilon_N
\end{equation}
for all $N$. This gives us that
$$
\sqrt{\widehat{V}_j(\mbf{n}_N)}\widehat{P}_j(\mbf{n}_N) \geq \frac{1}{\sqrt{2\pi}} \exp\left(-\frac{\widehat{M}_j(\mbf{n}_N)^2}{2 \widehat{V}_j(\mbf{n}_N)}\right) - \epsilon_N.
$$
Now note that
$$
\lim_{N\to\infty} \frac{1}{|\mbf{n}_N|}\widehat{V}_j(\mbf{n}_N)
=
\lim_{N\to\infty} \sum\limits_{k=1}^m \frac{n_k}{|\mbf{n}_N|} \Var[\widehat{X}'_{kj}] + \frac{1}{|\mbf{n}_N|}\Var[Y_{k, j}]
= 
\sum\limits_{k=1}^m \rho_k \Var[\widehat{X}'_{kj}] > 0.
$$
Furthermore, using that $\sum_{k=1}^m \rho_k\E[\widehat{X}'_{kj}] = 0$ and Assumption \ref{ass:sequence_assumption}, we have
$$
\left|\widehat{M}_j(\mbf{n}_N)\right| \leq \left|\sum\limits_{k=1}^m \left(n_k - |\mbf{n}_N|\rho_k\right) \E[\widehat{X}'_{kj}]\right| + \sum\limits_{k=1}^m|\E[Y_{k, j}]| \leq \tilde C\sqrt{|\mbf{n}_N|}.
$$
Combining the above two statements, we find that
$$
\frac{\widehat{M}_j(\mbf{n}_N)^2}{2 \widehat{V}_j(\mbf{n}_N)}
$$
is bounded, so that
\begin{equation}
\liminf_{N\to\infty} \frac{1}{\sqrt{2\pi}} \exp\left(-\frac{\widehat{M}_j(\mbf{n}_N)^2}{2 \widehat{V}_j(\mbf{n}_N)}\right) > 0.
\end{equation}
Since $\epsilon_N \to 0$, we thus find that 
$$
\frac{1}{\sqrt{2\pi}} \exp\left(-\frac{\widehat{M}_j(\mbf{n}_N)^2}{2 \widehat{V}_j(\mbf{n}_N)}\right) - \epsilon_N > \delta > 0
$$
for $N$ large enough. This then gives us that
\begin{equation}
    \liminf_{N\to\infty} \frac{1}{|\mbf{n}_N|}\log\left(\sqrt{\widehat{V}_j(\mbf{n}_N)}\widehat{P}_j(\mbf{n}_N)\right) \geq  \liminf_{N\to\infty} \frac{1}{|\mbf{n}_N|}\log\delta = 0.
\end{equation}
Note furthermore that
\begin{align*}
    \liminf_{N\to\infty} \frac{1}{|\mbf{n}_N|}\log\left(\sqrt{\widehat{V}_j(\mbf{n}_N)}\widehat{P}_j(\mbf{n}_N)\right)
    &=
    \liminf_{N\to\infty} \frac{1}{2|\mbf{n}_N|}\log\widehat{V}_j(\mbf{n}_N) + \frac{1}{|\mbf{n}_N|}\log\widehat{P}_j(\mbf{n}_N)
    \\
    &=
    \liminf_{N\to\infty} \frac{1}{|\mbf{n}_N|}\log\widehat{P}_j(\mbf{n}_N).
\end{align*}
Here we used that 
$$
\lim_{N\to\infty} \frac{1}{|\mbf{n}_N|}\widehat{V}_j(\mbf{n}_N) > 0,
$$
from which it follows that
$$
\lim_{N\to\infty} \frac{1}{|\mbf{n}_N|}\log\widehat{V}_j(\mbf{n}_N) = 0.
$$
Collecting everything, we obtain that
$$
\liminf_{N\to\infty} \frac{1}{|\mbf{n}_N|}\log\widehat{P}_j(\mbf{n}_N) \geq 0
$$
as desired. Finally, by using \eqref{eq:lower_bound_2} it follows that 
\begin{equation}
    \begin{split}
        \liminf_{N \to \infty}\frac{1}{|\mbf{n}_N|} \log P(\mbf{n}_N) &\geq \sum\limits_{j, k = 1}^m \rho_k \log \E[e^{\tau_j X'_{kj}}]= \sum\limits_{k=1}^m \rho_k \log \E[e^{\boldsymbol{\tau} \cdot \mbf{X}'_k}].\\
    \end{split}
\end{equation}
Recalling that $\boldsymbol{\tau}$ is the the unique vector of minimizers of $(\phi_1(\lambda_1), \ldots, \phi_m(\lambda_m))$, we see that the lower bound is proven, since
\begin{equation}
    \begin{split}
        \sum\limits_{k=1}^m \rho_k \log \E[e^{\boldsymbol{\tau} \cdot \mathbf{X}_k'}] &= \sum\limits_{k=1}^m \rho_k \log \E[e^{\boldsymbol{\tau} \cdot (\mathbf{X}_k - \mathbf{e}_k)}]\\
        &= \sum\limits_{k=1}^m - \tau_k \rho_k + \log\mathbb{E}[e^{\boldsymbol{\tau} \cdot \mathbf{X}_k}]\\
        &= - \sup_{\boldsymbol{\lambda}} \{\boldsymbol{\lambda} \cdot \boldsymbol{\rho} - \sum\limits_{k=1}^m \rho_k \log \E[e^{\boldsymbol{\lambda} \cdot \mathbf{X}_k}]\}.
    \end{split}
\end{equation}
\end{proof}

\section*{Acknowledgements}
This publication is part of the project ``Random graph representation of nonlinear evolution problems"  of the research programme Mathematics Cluster/NDNS+ which is partly financed by the Netherlands Research Organisation (NWO). IK also gratefully acknowledges support form NWO research program VIDI, project number VI.Vidi.213.108.

\bibliographystyle{alpha}
\bibliography{references.bib}

\end{document}